\documentclass[a4paper,reqno]{amsart}

\usepackage{amsfonts}

\title[Gradient bounds for minimizers of free discontinuity problems]{Gradient bounds for minimizers\\ of free discontinuity problems related to\\ cohesive zone models in fracture mechanics}

\author{G.~Dal Maso}
\author{A.~Garroni}

\address[Gianni Dal Maso]{SISSA
\\via Beirut 2-4\\
34014 Trieste\\
Italy}
\email[Gianni Dal Maso]{dalmaso@sissa.it}

\address[Adriana Garroni]{Dipartimento di 
Matematica\\
Universit\`a di Roma ``La Sapienza''\\
\hbox{  \ \ \ \ }
 \hbox{  \ \ \ \  \ \ } P.le Aldo Moro 3
\\00185 Roma\\
Italy}
\email{garroni@mat.uniroma1.it}

\newcommand{\lbl}[1]{\label{#1}}                

\newtheorem{thm}{Theorem}                    
\newcommand{\bthm}{\begin{thm}}
\newcommand{\ethm}{\end{thm}}

\newcommand{\nl}{\mbox{}\\}                     
\newcommand{\nll}{\nl\nl}                       

\newcommand{\beq}{\begin{eqnarray*}}                    
\newcommand{\eeq}{\end{eqnarray*}}

\newcommand{\beqn}[1]{\begin{eqnarray}\lbl{#1}}         
\newcommand{\eeqn}{\end{eqnarray}}

\newcommand{\barr}[1]{\begin{array}{#1}}        
\newcommand{\earr}{\end{array}}

\newcommand{\eps}{\varepsilon}
\newcommand{\e}{\varepsilon}                   
\newcommand{\Om}{\Omega}                        
\newcommand{\R}{{\mathbb R}}
\newcommand{\rn}{{\R^n}}       
\newcommand{\E}{{\mathcal E}}   
\newcommand{\into}[1]{\int_\Om #1}             
\newcommand{\hau}{{\mathcal H}^{n-1}}

\begin{document}

\begin{abstract}
\noindent
In this note we consider a free discontinuity problem for a scalar function, whose energy depends also on the size of the jump. We prove that the gradient of every smooth local minimizer never exceeds a constant, determined only by the data of the problem.
\nll
Keywords: fracture mechanics, cohesive zone models, functions of bounded variations, local minimizers, free discontinuity problems.
\nll
Mathematical Subjects Classification: 
74R99
, 
49J45
,
49N60
.
\end{abstract}

\maketitle


\section{Introduction}
\lbl{s-intro}

The study of cohesive zone models in fracture mechanics in the one dimensional case (see, e.g., \cite{Tru} and \cite{DelP})  leads to functionals of the form 
\beqn{f-oned}
\int_0^l F(|\dot u|)\, dx + \sum_{S(u)}G(|[u]|)\qquad u\in SBV(0,l)\,,
\eeqn
where $F\colon [0,+\infty)\to [0,+\infty)$ is $C^1$, strictly convex, increasing, superlinear at infinity, and satisfies $F(0)=F'(0)=0$, and $G\colon [0,+\infty)\to  [0,+\infty)$ is $C^1$, concave, and satisfies $G(0)=0$ and $G'(0)>0$. Here and in the rest of the paper $SBV$ is the space of special functions with bounded variation, for which we refer to \cite{AmFuPa}, $S(u)$ denotes the jump set of $u$, and $[u]$ denotes the jump of $u$.

To prove the existence of a minimizer of (\ref{f-oned}) with  appropriate boundary conditions we can consider the corresponding relaxed functional in $L^1(0,l)$, which for every $u\in BV(0,l)$ can be written as
\beqn{f-onedrelax}
\int_0^l \overline F(|\dot u|)\, dx + \sum_{S(u)}G(|[u]|) + G'(0)\, |u'_c|(0,l) \,,
\eeqn
where $\dot u$ is the density of the absolutely continuous part of the distributional derivative $u'$ and $u'_c$  is its Cantor part. In (\ref{f-onedrelax})  $\overline F(\xi)=F(\xi)$ for $\xi\leq e_M$ and $\overline F(\xi) = F(e_M) + F'(e_M)\,(\xi-e_M)$ if $\xi>e_M$, where $e_M$ is the unique constant such that $F'(e_M)=G'(0)$.
It is possible to prove that the minimum problem for the relaxed functional (\ref{f-onedrelax}) with appropriate boundary conditions has a solution. Moreover in \cite{BraDMGar} it was proved, by using one dimensional arguments, that if $G$ is strictly concave, then every local minimizer $u$ of (\ref{f-onedrelax}) satisfies 
$$
|\dot u|\leq e_M\quad \hbox{a.e. on }\ (0,l)\ ,\qquad |u'_c|(0,l)=0\,. 
$$
In particular this implies that $\overline F(\dot u)=F(\dot u)$ a.e. on $(0,l)$, so that $u$ is a local minimizer of~(\ref{f-oned}). Moreover
$$
F'(|\dot u|) \leq G'(0)\,.
$$
This justifies the interpretation of $G'(0)$ as the ultimate stress for the problem (see, e.g., \cite{Car}).

In this note we study the same problem in dimension $n\ge 1$.
We consider functionals of the form
\beqn{f1}
\into F(|\nabla u|) \,dx  + \int_{S(u)} G(|[u]|) \,
d\hau\qquad u\in SBV(\Om)\,,
\eeqn
where $\nabla u$ is the density of the absolutely continuous part of the distributional gradient $Du$, and $F$ and $G$ satisfy the same properties considered for (\ref{f-oned}).

Also in this case the functional is not lower semicontinuous, so in order to prove existence results we consider its {\it relaxed functional\/} in $L^1(\Om)$ (see~\cite{BouBrBu}), which  is represented on $BV(\Om)$ by 
\beqn{f1relax}
\E(u)=\into \overline F(|\nabla u|) \,dx  + \int_{S(u)} G(|[u]|) \,
d\hau + G'(0)\, |D^cu|(\Om)\,,
\eeqn
where $\overline F$ is defined as for (\ref{f-onedrelax}) and $D^cu$ denote the Cantor part of $Du$. Under appropriate boundary conditions the minimum problems for (\ref{f1relax}) have a solution. A {\it local minimizer\/} of ${\E}$ in $\Om$ is a function $u\in BV(\Om)$, with ${\E}(u)<+\infty$,  for which there exists $\eta>0$ such that $\E(u)\leq \E(v)$ for every $v\in BV(\Om)$ with ${\rm supp}(v-u)\subset\subset\Om$ and $\|v-u\|_{BV(\Om)}<\eta$.

Also in this case it is reasonable to  expect that  any local minimizer $u$ satisfies
\beqn{f-optimalcond2}
|\nabla u|\leq e_M\quad \hbox{a.e.\ on }\ \Om\ ,\qquad |D^cu|(\Om)=0\,,
\eeqn
where $e_M$ is defined as for (\ref{f-onedrelax}). In fracture mechanics the functionals (\ref{f1}) and (\ref{f1relax}) are used to study cohesive zone models in the antiplane case. In this context the first inequality in (\ref{f-optimalcond2}) says that the norm of the deformation gradient of the solution cannot exceed the constant  $e_M$, which is interpreted as the yield strain of the problem. Since (\ref{f-optimalcond2}) implies $F'(|\nabla u|)\leq G'(0)$ a.e.\ on $\Om$, the constant $G'(0)$ plays the role of the ultimate stress for the crack problem.

The aim of this note is to present a partial result in this direction. Namely, we prove that, if 
$$
\lim_{t\to 0+}\frac{G(t)-G'(0)\, t}{t^2}<0
$$
and $u$ is a local minimizer of (\ref{f1relax}) in $\Om$, then 
$$
|\nabla u|\leq e_M
$$
in every open subset of $\Om$ where $u$ is of class~$C^1$.
As a consequence we have that, if $u$ is a $C^1$ local minimizer for (\ref{f1relax}) in $\Om$, then it is also a local minimizer for (\ref{f1}).

\section{Statement and proof of the result}
\lbl{s-result}

Let $\Om$ be an open subset of $\rn$, $n\geq 1$. We assume that the functions $F$ and $G$ satisfy the following properties:
\begin{itemize}
\item[(a)] $F$ is $C^1$, strictly convex, increasing, and superlinear at infinity, and satisfies $F(0)=F'(0)=0$;
\item[(b)] $G$ is $C^{1}$, nonnegative, concave, and satisfies $G(0)=0$, $G'(0)>0$, and
\beqn{secder}
\lim_{t\to 0+}\frac{G(t)-G'(0)\, t}{t^2}<0\,.
\eeqn
\end{itemize}
The function $\overline F$ is defined as follows
\beqn{fbar}
\overline F(\xi) = \begin{cases}
F(\xi) \ \ \ & \text{if }\xi\leq e_M\,,\\
F(e_M) +  F'(e_M)\, (\xi-e_M) & \text{if }\xi > e_M
\,,
\end{cases}
\eeqn
where $e_M$ is the unique solution of the equation $F'(e_M)=G'(0)$.
\bigskip

\bthm
\lbl{theorem}
Assume that $F$ and $G$ satisfy conditions (a) and (b) and let $u$ be a local minimizer of the functional ${\E}$ defined by~(\ref{f1relax}). Suppose that $u$ is of class $C^1$ on an open subset  $U$ of $\Omega$. Then $|\nabla u|\leq e_M$  in $U$.
\ethm

The result stated in Theorem~\ref{theorem} implies that, if $u$ is a local minimizer of (\ref{f1relax}) satisfying $u\in C^1(\Om\setminus K)$, with $K$ closed and $\hau(K)<+\infty$, then $u$ is also a local minimizer of (\ref{f1}). Indeed in this case $D^cu=0$, hence $u\in SBV(\Om)$, and $\bar F(|\nabla u|) = F(|\nabla u|)$ a.e.\ in $\Om$ by Theorem~\ref{theorem}.

\begin{proof}[Proof of Theorem~\ref{theorem}]
Without loss of generality we consider only the case $e_M=F'(e_M)=G'(0)=1$ and $U=\Om$. We argue by contradiction and we assume that there exists a point $x_0\in\Om$ such that $|\nabla u(x_0)|=\lambda$, with $\lambda>1$. By changing the coordinate system, it is not restrictive to assume that $x_0=0$, $u(0)=0$, and $\nabla u(0)=\lambda\, e_n$, where $e_n:=(0,\dots,0,1)$ is the last vector of the canonical basis of $\R^n$.

We want to construct a competitor $w$ by modifying $u$ in a small set $V\subset\subset \Om$ with piecewise $C^1$ boundary in such a way that $w$ is close to $u$ in the $BV$ norm and the energy of $w$ is strictly  below the energy of $u$, contradicting the local minimality. In all cases we will take $w$ of the form
\beqn{a1}
w =\begin{cases}\alpha\, u\quad & \text{in }V\,,\\
u & \text{otherwise,}
\end{cases}
\eeqn
for a suitable constant $\alpha<1$. The problem is reduced to choose $\alpha$ and $V$ such that
\beqn{contraddizione}
\|u-w\|_{BV(\Om)}<\eta \qquad \hbox{and}\qquad
{\E}(u)-{\E}(w)>0\,,
\eeqn
where $\eta$ is the constant in the definition of local minimality for~$u$.

We consider three cases corresponding to different hypotheses on $G$ and $u$ with increasing level of difficulty.
\medskip

\noindent{\it Case 1: \ $G''(0)=-\infty$.}

\noindent Let us first consider the case where $G$ satisfies the following condition
\beqn{menoinfinito}
\lim_{t\to 0+}\frac{G(t)-t}{t^2}=-\infty\,.
\eeqn

Let us fix $\eps\in(0,\frac{1}{2})$, with $\lambda-\eps>1$. By the continuity of $\nabla u$ we can find $R>0$ small enough so that 
\beqn{star}
|\nabla u-\lambda\,e_n|<\eps\quad\hbox{in }B_R\,,
\eeqn
where $B_R$ is the closed ball with center~$0$ and radius~$R$. As a consequence we can show that $|\nabla u|>\lambda-\eps$ in $B_R$ and that there exists $\delta>0$ such that
\begin{eqnarray*}
& u(x)>\delta\quad \text{for every }\, x\in B_R\, \text{ with }\, x_n=\e R\,,
\\
& u(x)<-\delta\quad \text{for every }\, x\in B_R\, \text{ with }\, x_n=-\e R\,.
\end{eqnarray*}
This implies that for $0< \sigma<\delta$ the projection of the set $\{x\in B_R\,:\ u(x)=\sigma\}$ onto the hyperplane $\{x_n=0\}$ contains the projection of the set 
$\{x\in B_R\,:\ x_n=\e R\}$, and therefore
\beqn{*}
\hau({B_R\cap \{u=\sigma\}})\geq K_{\eps,R}:= \omega_{n-1}R^{n-1}({1-\eps^2})^{(n-1)/2}\,,
\eeqn
where $\omega_{n-1}$ is the $(n-1)$-dimensional measure of the unit ball in $\R^{n-1}$.
Moreover (\ref{star}) implies that there exists a constant  $L<+\infty$  such that
\beqn{disuguaglianza}
\hau(\{x\in \partial  B_R: 0<u(x)<\sigma\})\leq L\,\sigma
\eeqn
for every $\sigma> 0$.

For $0<\sigma<\delta$ we define
$$
V_\sigma:=\{x\in  B_R: 0<u(x)<\sigma\}\,.
$$
Since $u$ is $C^1$, there exists a constant $M$ such that 
$$
\hau(\partial V_\sigma)\le M
$$
for $0<\sigma<\delta$.

We now fix $\alpha<1$ such that $\alpha\,(\lambda-\eps)>1$ and 
$(1-\alpha)\,(\|u\|_{BV(\Om)}+\delta\,M)<\eta$,
and define $w$ as in (\ref{a1}) with $V:=V_\sigma$ for some $\sigma\in (0,\delta)$ to be chosen later. 
Since 
$$
\|w-u\|_{BV(\Om)}\le (1-\alpha)\,\|u\|_{BV(\Om)}+(1-\alpha)\,\sigma\,\hau(\partial V_\sigma)\,,
$$
we have $\|w-u\|_{BV(\Om)}\le\eta$ for $0<\sigma<\delta$, so that the first inequality in (\ref{contraddizione}) is satisfied.

Using the definition of $\mathcal E$ and $\bar F$, we get
\beqn{primocaso}
\nonumber{\E}(u)-{\E}(w) \!\!\! &=&\!\!\!  (1-\alpha)\int_{V_\sigma}|\nabla u|\,dx
-\int_{ B_R\cap \{u=\sigma\}} G((1-\alpha)u)\, d\hau
\\
&& \!\!\! {}
 - \int_{\partial B_R\cap V_\sigma} G((1-\alpha)u)\, d\hau\,.
\eeqn
Since $u$ is a $C^1$ local minimum of ${\E}$ and $|\nabla u|>1$ in $B_R$, in particular $u$ is a $C^1$ local
minimum of
$$
\int_{B_R}|\nabla u| \,dx
$$
and then, it satisfies  the Euler equation 
\beqn{div0}
{\rm div}\,
\Big(\frac{\nabla u}{|\nabla u|}\Big)=0
\quad\hbox{in the sense of distributions on }B_R \,.
\eeqn
Thus, by the divergence theorem, we have
\beqn{perparti}
\nonumber\int_{V_\sigma}|\nabla u|\,
dx \!\!\! &=&\!\!\!    \int_{B_R\cap \{u=\sigma\}} u \,d\hau +\int_{\partial B_R\cap V_\sigma}\frac{\nabla u}{|\nabla u|}\frac{x}{|x|} u
\,d\hau \\
\!\!\! &\geq&\!\!\!   \int_{B_R\cap \{u=\sigma\}} u \,d\hau -\int_{\partial B_R\cap V_\sigma} u
\,d\hau
\,.
\eeqn
Moreover, by condition (\ref{menoinfinito}), for any given $c>0$ we can choose $\sigma$ small enough so that
$$
G((1-\alpha)u)< (1-\alpha)u-c\,(1-\alpha)^2u^2\quad\text{on }\,V_\sigma\,.
$$
This, together with (\ref{perparti}) and (\ref{primocaso}), implies
\beqn{fineprimocaso}
\nonumber{\E}(u)-{\E}(w)\!\!\! &\geq&\!\!\! (1-\alpha)\sigma\,\hau(B_R\cap \{u=\sigma\})
-(1-\alpha)\int_{\partial B_R\cap V_\sigma}u
\,d\hau\\
\nonumber\  &&\!\! \!{}-(1-\alpha)\, \hau({B_R\cap \{u=\sigma\}})\, [\sigma-c\,(1-\alpha)\sigma^2]\\
\nonumber\ &&\!\!\! {}-
(1-\alpha)\int_{\partial B_R\cap V_\sigma}[u-c\,(1-\alpha)\,u^2] \,d\hau
\\
\nonumber&\geq&\!\!\! (1-\alpha)\,\sigma\,[c\,(1-\alpha)\,\sigma\,\hau({B_R\cap \{u=\sigma\}})- 2\,\hau(\partial B_R\cap V_\sigma)]\,.
\eeqn
{}From (\ref{*}) and (\ref{disuguaglianza}) we get
$$
{\E}(u)-{\E}(w)\geq (1-\alpha)\,\sigma^2[c\,(1-\alpha)\,K_{\eps,R}-2L]\,,
$$
which gives the second inequality in (\ref{contraddizione}) when $c$ is big enough.
\medskip

Next we consider the general case where $G$ does not necessarily satisfy (\ref{menoinfinito}). In this case we must choose the set $V$ more carefully. 
In order to explain the new ideas of the proof without technicalities, we prove first  the result in two dimensions in the simplest case: when $u$ is an affine function.
\medskip

\goodbreak
\noindent{\it Case 2:\ \ $-\infty<G''(0)<0$, $u$ affine, and $n=2$.}
\par
\nopagebreak
\noindent  We now consider the case $n=2$ with $u$ affine. We assume that $G$ satisfies the following condition
\beqn{secderfin}
-\infty<\lim_{t\to 0+}\frac{G(t)-t}{t^2}<0\,.
\eeqn
Then there exist two constants $c_2>c_1>0$ such that
\beqn{secderfin*}
t-c_2\, t^2< G(t)<t- c_1 t^2
\eeqn
for $t>0$ small enough.

It is not restrictive to take $u(x)=\lambda x_2$ for every $x=(x_1,x_2)\in\Om\subseteq\R^2$. We assume by contradiction that  $\lambda>1$.
 It is easy to check that in general we may not choose $V$ to be a rectangle. Indeed, if $V=\{(x_1,x_2)\in\Om\,:\ 0<x_1<S\,,\ 0<x_2<\delta\}$, following the computation of Case~1  we get  for $\delta>0$ small enough
\beqn{rettangoloaffine}
\nonumber \E(u)-\E(w)\!\!\!&\leq&\!\!\!(1-\alpha)\lambda S\delta-
\int_{0}^{\delta}[(1-\alpha)\lambda x_2 - c_2(1-\alpha)^2\lambda^2x^2_2]\,dx_2\\
\nonumber &&\!\!\!{}- (1-\alpha)\lambda \delta S +c_2 S(1-\alpha)^2\lambda^2\delta^2\\
\nonumber &=&\!\!\!{}-(1-\alpha)\lambda\frac{\delta^2}{2} +c_2(1-\alpha)^2\lambda^2\frac{\delta^3}{3}+c_2S(1-\alpha)^2\lambda^2\delta^2\,,
\eeqn
and the right-hand side is positive for every $\delta>0$ only if $S\ge [2\,(1-\alpha)\,\lambda\, c_2]^{-1}$. This condition may be incompatible with the inclusion $V\subset\subset\Om$. For the same reason we can not define $V$ as in Case~1.

Since the previous computation shows that the problem is given by the short sides of the rectangle, we are led to overcome this difficulty by defining a special profile for the boundary of $V$. Let us fix $r$ and $R$, with $r<R$, and let $\varphi\colon [0,R]\to [0,+\infty)$ be a nonincreasing function, to be chosen later, satisfying $\varphi(\rho)=1$ in $0\leq \rho\leq r$ and $\varphi(R)=0$. We take $V$ of the form
$$
V:=\{(x_1,x_2)\ :\ |x_1|<R\,,\ 0<x_2<\sigma\varphi(|x_1|)\}\,,
$$
with $0<\sigma<1$, and we consider the function $w$ defined by~(\ref{a1}).
Let us compute the energy of $w$ and show that (\ref{contraddizione}) holds for a suitable choice of $r$, $R$, $\varphi$, $\sigma$, and $\alpha$.

If $\alpha<1$ and $\alpha\,\lambda>1$,
using the definition of $w$ we get
\beqn{e0}
 \nonumber{\E}(u)-{\E}(w)\!\!\!&=&\!\!\!(1-\alpha)\,\lambda\,{\mathcal L}^2(V) -\int_{\partial V\setminus\{x_2=0\}} G((1-\alpha)\,\lambda\,x_2) d{\mathcal H}^1(x)
 \\
 \nonumber &=&\!\!\! 2\,(1-\alpha)\,\lambda\,r\,\sigma+2\,(1-\alpha)\,\lambda\int_r^R\sigma\,\varphi(\rho)\,d\rho 
 -2\,r \,G((1-\alpha)\,\lambda\,\sigma)
 \\
 \nonumber&&\!\!\!{}-2\int_r^R G((1-\alpha)\,\lambda\,\sigma\,\varphi(\rho))\sqrt{1+(\sigma\,\varphi'(\rho))^2}\,d\rho\,.
 \eeqn
Using the fact that $\sqrt{1+t^2}\leq 1+\frac{1}{2} t^2$ and $0\leq G(t)\leq t-c_1t^2$ for small $t>0$ we obtain 
\beqn{e00}
 \nonumber
\frac{1}{(1-\alpha)\lambda}\bigl({\E}(u)-{\E}(w)\bigr)\!\!\!&\geq&\!\!\! 2\,c_1r\,(1-\alpha)\,\lambda \,\sigma^2
-\int_r^R\sigma^3\varphi(\rho)\,(\varphi'(\rho))^2\,d\rho \\
 \nonumber&&\!\!\!{}+2\,c_1\int_r^R(1-\alpha)\,\lambda\,\sigma^2\,(\varphi(\rho))^2\, d\rho
 \\
 \nonumber&\geq&\!\!\! \int_r^R[2\,c_1(1-\alpha)\,\lambda\,\sigma^2 (\varphi(\rho))^2-\sigma^3\,\varphi(\rho)\, (\varphi'(\rho))^2]\,d\rho\,.
 \eeqn
The inequality ${\E}(u)-{\E}(w)>0$ can be obtained easily for $\sigma$ small enough if
$$
\varphi'(\rho)^2\varphi(\rho)=k\,(\varphi(\rho))^2
$$
for a suitable constant $k$ independent of $\sigma$. It is easy to check that a solution of this equation on $[r,R]$, with $\varphi(r)=1$ and $\varphi(R)=0$, is given by
\beqn{profilo}
\varphi(\rho)=\frac{(\rho-R)^2}{(r-R)^2}\,,
\eeqn
with $k=4\,(R-r)^{-2}$.
With this choice of the profile $\varphi$ we get 
\beqn{stimadalbasso}
{\E}(u)-{\E}(w)\geq (1-\alpha)\,\lambda \int_r^R
 [2\,c_1(1-\alpha)\,\lambda\,\sigma^2- 4\,\sigma^3(R-r)^{-2}]\,(\varphi(\rho))^2 d\rho\,.
\eeqn
Now we choose $\alpha<1$ such that $\alpha\,\lambda>1$ and
$$
(1-\alpha)\,\Big[\|u\|_{BV(\Om)}+ 2\,r\,\lambda + 2\,\lambda \int_r^R \varphi(\rho))\sqrt{1+(\varphi'(\rho))^2}\,d\rho \Big] < \eta\,.
$$
Since
$$
\|w-u\|_{BV(\Om)} \le (1-\alpha)\,\Big[\|u\|_{BV(\Om)}+ 2\,\sigma\, r\,\lambda + 2 \, \lambda \int_r^R \sigma \,\varphi(\rho))\sqrt{1+(\sigma\,\varphi'(\rho))^2}\,d\rho \Big]\,,
$$
the first inequality in (\ref{contraddizione}) is satisfied for $0<\sigma<1$. By (\ref{stimadalbasso}) the second  inequality in (\ref{contraddizione}) is satisfied for $0<\sigma< c_1(1-\alpha)\,\lambda\,(R-r)^2/2$.
This concludes the proof of Case~2.
 \medskip
 
\noindent{\it Case 3: General case.}
 
\noindent We finally prove the result in the general case.  As in Case~1, for a given $\e\in(0,\frac{1}{2})$
such that $\lambda-\e>1$ we may select $R>0$ so small that $|\nabla u-\lambda\,e_n|<\eps$ and $|\nabla u|>\lambda-\eps$ in $B_R$.
Now, inspired by the calculation of Case~2, we fix  $r>0$, with $r<R$, and we consider
the function $a(x)$ defined in $B_R$ by $a(x)=\varphi(|x|)$; i.e.,
$$
a(x)=\begin{cases}
\displaystyle{\frac{(|x|-R)^2}{(r-R)^2}}\quad &\text{if }
r<|x|<R\,,
\vspace{5pt}
\\
1 & \text{if }|x|\leq r\,.
\end{cases}
$$
Let $v:={u/a}$ and $S_\sigma =:\{x\in B_R: v(x)=\sigma\}=\{x\in B_R: u(x)=\sigma\,a(x)\}$. Since $u$ is $C^1$, there exist $\delta>0$ and $M>0$ such that
\beqn{e101}
\hau(S_\sigma)\le M
\eeqn
for $0<\sigma<\delta$.

We now fix $\alpha<1$ such that $\alpha\,(\lambda-\eps)>1$ and 
$(1-\alpha)\,[\|u\|_{BV(\Om)}+\delta\,M]<\eta$,
and define $w$ as in (\ref{a1}) with $V:=\{{x\in B_R}:\ 0\leq v(x)\leq\sigma\}$.
Since 
$$
\|w-u\|_{BV(\Om)}\le (1-\alpha)\,\|u\|_{BV(\Om)}+
(1-\alpha)\,\sigma\,\hau(S_\sigma)\,,
$$
we have $\|w-u\|_{BV(\Om)}\le\eta$ for $0<\sigma<\delta$, so that the first inequality in (\ref{contraddizione}) is satisfied.

To conclude the  proof we have to show that $\sigma$ can be chosen in $(0,\delta)$ so  that the second inequality in (\ref{contraddizione}) holds, contradicting the local minimality of~$u$. If $\delta$ is small enough, we may assume that $G$ satisfies the second inequality of (\ref{secderfin*}) for $0<t<\delta$.
Let
$C_R^r:=B_R\setminus B_r$. By the definition
of $w$ we have $|\nabla w|= \alpha |\nabla u|>1$ a.e.\ on $V$ and thus 
\beqn{e5}
\nonumber &&{\E}(u)-{\E}(w) = (1-\alpha)\int_{V}|\nabla u|\,
dx - \int_{B_R\cap S_\sigma} G((1-\alpha)u)\, d\hau
\\
\nonumber &\geq&(1-\alpha)\int_{V}|\nabla u|\,
dx - (1-\alpha)\int_{B_r\cap S_\sigma} u\, d\hau -
(1-\alpha)\int_{C_R^r\cap S_\sigma} u\, d\hau 
\\
&& {}+ c_1(1-\alpha)^2\int_{C_R^r\cap S_\sigma} u^2\, d\hau+c_1(1-\alpha)^2\int_{ B_r\cap S_\sigma} u^2\, d\hau\,.
\eeqn

As in Case 1 we use the fact that $u$ satisfies (\ref{div0}). Since $\frac{\nabla v}{|\nabla v|}$ is the outer unit normal to $S_\sigma$ and $\nabla v=\nabla u$ on $B_r\cap S_\sigma$,
by the divergence theorem we get
\beqn{e6}
\int_{V}|\nabla u|\,
dx =  \int_{C_R^r\cap S_\sigma}\frac{\nabla u}{|\nabla u|} \frac{\nabla v}{|\nabla v|} u
\,d\hau + \int_{B_r\cap S_\sigma} u \,d\hau\,.
\eeqn
Since $\nabla v={\nabla u/
a} -{u\nabla a/ a^2}=(1/a)(\nabla u-\sigma \nabla a)$ on $C_R^r\cap S_\sigma$, we have
$$
\frac{\nabla u}{|\nabla u|} \frac{\nabla v}{|\nabla v|} =\frac{ \displaystyle |\nabla u| - \sigma \frac{\nabla a\cdot\nabla u }{ |\nabla u|} }{  \displaystyle
|\nabla u -\sigma \nabla a|}
$$
on $C_R^r\cap S_\sigma$.
Using Taylor's expansion of the right-hand side with respect to $\sigma$ we obtain 
$$
\frac{\nabla u}{|\nabla u|} \frac{\nabla v}{|\nabla v|} = 1 + 
 \frac{(\nabla a\cdot \nabla u)^2-|\nabla a|^2 | \nabla u|^2 }{ 2\, |\nabla u|^4}\sigma^2 + O(\sigma^3)\,,
$$
and hence
\beqn{e7}
\frac{\nabla u}{|\nabla u|} \frac{\nabla v}{|\nabla v|} \geq 1 - \frac{|\nabla a|^2}{2\,(\lambda-\eps)^2} \sigma^2 + 
O(\sigma^3)
\eeqn
on $C_R^r\cap S_\sigma$.
Since $|\nabla a|^2 a=4\,(R-r)^{-2} a^2$ on $C_R^r\cap S_\sigma$, by (\ref{e5}), (\ref{e6}), and (\ref{e7}) we have
$$
{\E}(u)-{\E}(w) \geq \sigma^2 (1-\alpha)\,[c_1(1-\alpha) - K_{\eps,r,R}\,\sigma] 
\int_{C_R^r\cap S_\sigma} a^2\, d\hau + O(\sigma^4) \,,
$$
with $K_{\eps,r,R}:= 2\,(R-r)^{-2}(\lambda-\eps)^{-2}$. 
Taking now $\sigma>0$ small enough we obtain ${\E}(u)-{\E}(w)>0$, which concludes the proof.
\end{proof}

{\frenchspacing
\normalfont 
}
\end{document}